\documentclass[11pt]{article}

\usepackage[margin=1.25in]{geometry}

\usepackage{amssymb,amsmath}
\usepackage{amsthm}
\usepackage{hyperref}
\usepackage{mathrsfs}
\usepackage{textcomp}
\hypersetup{
    colorlinks,%
    citecolor=black,%
    filecolor=black,%
    linkcolor=black,%
    urlcolor=black
}

\usepackage{tikz}
\usetikzlibrary{matrix,arrows}

\usepackage{verbatim}

\usepackage{sectsty}

\makeatletter

\newdimen\bibspace
\renewenvironment{thebibliography}[1]{%
 \section*{\refname 
       \@mkboth{\MakeUppercase\refname}{\MakeUppercase\refname}}%
     \list{\@biblabel{\@arabic\c@enumiv}}%
          {\settowidth\labelwidth{\@biblabel{#1}}%
           \leftmargin\labelwidth
           \advance\leftmargin\labelsep
           \itemsep\bibspace
           \parsep\z@skip     %
           \@openbib@code
           \usecounter{enumiv}%
           \let\p@enumiv\@empty
           \renewcommand\theenumiv{\@arabic\c@enumiv}}%
     \sloppy\clubpenalty4000\widowpenalty4000%
     \sfcode`\.\@m}
    {\def\@noitemerr
      {\@latex@warning{Empty `thebibliography' environment}}%
     \endlist}

\makeatother

\makeatletter

\newtheorem{thm}{Theorem}[section]
\newtheorem{lem}[thm]{Lemma}
\newtheorem{prop}[thm]{Proposition}

\newtheorem{conj}[thm]{Conjecture}
\newtheorem{cor}[thm]{Corollary}
\newtheorem{rem}[thm]{Remark}

\def\XXint#1#2#3{{\setbox0=\hbox{$#1{#2#3}{\int}$}
  \vcenter{\hbox{$#2#3$}}\kern-.5\wd0}}

                \newcommand{\lda}{\lambda}
\newcommand{\om}{\Omega}                \newcommand{\pa}{\partial}
\newcommand{\va}{\varepsilon}           \newcommand{\ud}{\mathrm{d}}
\newcommand{\be}{\begin{equation}}      \newcommand{\ee}{\end{equation}}
\newcommand{\w}{\omega}                 
\newcommand{\Lda}{\Lambda}              
\newcommand{\R}{\mathbb{R}}

\begin{document}

\title{\textbf{On the isoperimetric quotient over scalar-flat\\ conformal classes}
\bigskip}

\author{\medskip  Tianling Jin\footnote{T. Jin is partially supported by Hong Kong RGC grant ECS 26300716 and HKUST initiation grant IGN16SC04.}, \  \
Jingang Xiong\footnote{J. Xiong is partially supported by NSFC 11501034, a key project of NSFC 11631002 and NSFC 11571019.}}

\date{\today}

\maketitle

\begin{abstract}
Let $(M,g)$ be a smooth compact Riemannian manifold of dimension $n$ with smooth boundary $\partial M$. Suppose that $(M,g)$ admits a scalar-flat conformal metric.   We prove that the supremum of the isoperimetric quotient  over the scalar-flat conformal class is strictly larger than the best constant of the isoperimetric inequality in the Euclidean space, and consequently is achieved, if either (i) $n\ge 12$ and $\partial  M$ has a nonumbilic point; or (ii) $n\ge 10$, $\partial M$ is umbilic and the Weyl tensor does not vanish at some boundary point. 
\end{abstract}

\section{Introduction}

In 1921, Carleman \cite{Carleman} proved that the classical isoperimetric inequality holds for a simply connected domain on a minimal surface,  by showing  the sharp inequality:
\begin{equation}\label{eq:carleman}
\int_B|f|^2\le\frac{1}{4\pi}\left(\int_{\partial B}|f|\right)^2
\end{equation}
 for every holomorphic function $f$ on the unit ball $B\subset\R^2$. Jacobs in \cite{J} extended \eqref{eq:carleman} to general bounded open set $\Omega\subset\R^2$ with smooth boundary $\partial\Omega$: there exists a positive constant $C_\Omega$ such that
\begin{equation}\label{eq:carleman2}
\int_\Omega|f|^2\le C_\Omega\left(\int_{\partial \Omega}|f|\right)^2
\end{equation}
 for every holomorphic function $f$ on $\Omega$.  Moreover, when $\Omega$ is not simply connected, the best constant $C_\Omega>\frac{1}{4\pi}$ and  is achieved. 
 
 A corollary of \eqref{eq:carleman} is the sharp inequality:
 \begin{equation}\label{eq:carleman0}
\int_B e^{2w}\le\frac{1}{4\pi}\left(\int_{\partial B}e^w\right)^2
\end{equation}
for every harmonic function $w$ on $B\subset\R^2$. 
Note that for a harmonic function $w$ on $B\subset\R^2$, the Gauss curvature of $(B, e^{w}g_{\R^2})$ is identically zero, where $g_{\R^2}$ is the Euclidean metric on $\R^2$.

In \cite{HWY1}, Hang-Wang-Yan obtained a higher dimensional generalization of the above inequalities on the unit ball $B\subset\R^n$ with $n\ge 3$. They proved the following sharp inequality:
\begin{equation}\label{eq:HWY}
\|v\|_{L^\frac{2n}{n-2}(B)}\le n^{-\frac{n-2}{2(n-1)}} \w_n^{-\frac{n-2}{2n(n-1)}}\|v\|_{L^\frac{2(n-1)}{n-2}(\partial B)}
\end{equation}
for every harmonic function $v$ on $B$ with $\omega_n$ being the volume of $B$, and classified all the minimizers. For a positive harmonic function $v$  on $B$, the scalar curvature of $(B, v^{\frac{4}{n-2}}g_{\R^n})$ is identically zero, where $g_{\R^n}$ is the Euclidean metric on $\R^n$. 

In \cite{HWY}, Hang-Wang-Yan further studied a generalization of \eqref{eq:HWY} on a smooth compact  Riemannian manifold $(M,g)$ of dimension $n\ge 3$ with smooth boundary $\pa M$, by considering the following variational problem on the isoperimetric quotient over the scalar-flat conformal class
  \be\label{eq:theproblem}
 \Theta_{M,g}=\sup_{\tilde g\in \mathcal{A}_g} \frac{\mbox{Volume}_{\tilde g} (M)}{\mbox{Area}_{\tilde g} (\pa M)^{\frac{n}{n-1}}},
 \ee
where  $\mathcal{A}_g:=\{\tilde g\in [g]:\mbox{ the scalar curvature }R_{\tilde g}=0\}$, and $[g]$ is the conformal class of $g$. It was explained in \cite{HWY} that the set $\mathcal{A}_g$ is not empty if and only if the first eigenvalue $\lda_1(L_g)$ of the conformal Laplacian operator
$
L_g:= -\Delta_g +\frac{n-2}{4(n-1)}R_g
$
with zero Dirichlet boundary condition is  positive. Note that the positivity of $\lda_1(L_g)$ does not depend on the choice of the metrics in $[g]$.
Assuming $\lda_1(L_g)>0$, they proved in \cite{HWY} that
  \begin{equation}\label{eq:highDisoperimetric}
\Theta_{\overline B_1,g_{\R^n}}\le \Theta_{M,g}<\infty,
\end{equation}
 and $\Theta_{\overline B_1,g_{\R^n} }$ coincides with the best constant of the isoperimetric inequality in the Euclidean space, that is,
 \[
\Theta_{\overline B_1,g_{\R^n} }= n^{-\frac{n}{n-1}} \w_n^{-\frac{1}{n-1}}.
\]
They also showed in \cite{HWY} that $\Theta_{M,g}$ is achieved if the strict inequality
\be \label{eq:mainstrict}
\Theta_{\overline B_1,g_{\R^n}}< \Theta_{M,g}
\ee
  holds, and made a conjecture that:

\begin{conj}[\cite{HWY}]\label{conjecture} Assume $n\ge 3$, $(M,g)$ is a  smooth compact Riemannian manifold of dimension $n$ with nonempty smooth boundary $\pa M$, and $\lda_1(L_g)>0$. If $(M,g)$ is not conformally diffeomorphic to $(\overline B_1, g_{\R^n})$, then the strict inequality  \eqref{eq:mainstrict} holds.

\end{conj}

In this paper, we prove

\begin{thm} \label{thm:main} Let $(M,g)$ be a  smooth compact Riemannian manifold of dimension $n$ with nonempty smooth boundary $\pa M$. Suppose that $\lda_1(L_g)>0$. Then the strict inequality \eqref{eq:mainstrict} holds if one of the following two conditions holds:
\begin{itemize}
\item[(i)] $n\ge 12$ and $\partial  M$ has a nonumbilic point;
\item[(ii)] $n\ge 10$, $\pa M$ is umbilic and the Weyl tensor $W_g\neq 0$ at some boundary point.
\end{itemize}
Therefore, $\Theta_{M,g}$ is achieved under one of these two conditions.
\end{thm}

Consequently, we have

\begin{cor}\label{cor:1} Let $\mathcal{O}$ be a bounded open subset of $\R^n$, $n\ge 12$, with a smooth connected boundary. Then $\Theta_{\overline B_1,g_{\R^n}}= \Theta_{\overline{\mathcal{O}},g_{\R^n}}$ if and only if $\mathcal{O}$ is a ball.
\end{cor}

At the time of writing this paper, we learned from Professor Meijun Zhu that, together with M. Gluck, they \cite{GZ} recently verified \eqref{eq:mainstrict}  when $M=\overline B_1\setminus  B_\va$ for sufficiently small $\va>0$ with the Euclidean metric in all dimensions.

Throughout the paper, we will always assume that $\lda_1(L_g)>0$. Denote the Poisson kernel of $L_g u=0$ with Dirichlet boundary condition by $P_{g}$. Then $P_{g}$ is  positive and satisfies the following conformal transformation law: 
\begin{equation}\label{eq:conformalchange}
P_{\rho^{\frac{4}{n-2}}g}v=\rho^{-1}P_{g}(\rho v), \quad \mbox{for all }v\in C^\infty(\partial M),\ \rho\in C^\infty(\overline M) \mbox{ and }\rho>0,
\end{equation}
where
\[
(P_{g} v)(x)=\int_{\partial M}P_{g}(x,\xi)v(\xi)\,\ud s_g(\xi).
\]
Making use of the conformal transformation law of the conformal Laplacian:
\[
L_{ \rho^{\frac{4}{n-2}}g} \phi= \rho^{-\frac{n+2}{n-2}} L_g (\rho\phi) \quad \mbox{for all }\rho, \phi \in C^\infty(\overline M) \mbox{ and }\rho>0,
\]
it was pointed out in \cite{HWY} that
\begin{equation}\label{eq:energyinf}
\Theta_{M,g}=\sup\left\{I[v]: v\in L^{\frac{2(n-1)}{n-2}}(\pa M), v\neq 0\right\},
\end{equation}
where
\[
I[v]=\frac{\int_{M} |P_{g} v|^{\frac{2n}{n-2}}\,\ud Vol_g}{(\int_{\pa M} |v|^{\frac{2(n-1)}{n-2}}\,\ud s_g)^{\frac{n}{n-1}}}.
\]

One sees from the definition that $\Theta_{M,g}$ depends only on the conformal class $[g]$. The results of the above variational problem \eqref{eq:theproblem} or \eqref{eq:energyinf} in \cite{HWY} and in this paper show an analogy to the Yamabe problem solved by Yamabe \cite{Yamabe}, Trudinger \cite{Trudger}, Aubin \cite{Aubin} and Schoen \cite{Schoen}, as well as to a boundary Yamabe problem (or higher dimensional Riemannian mapping problem) studied by Escobar \cite{E}, Marques \cite{M1, M2}, Almaraz \cite{A}, etc. Actually, we have the following interesting triangle diagram:

\centerline{\small
\begin{tikzpicture}[node distance=3cm, auto]
  \node (T) {$L^{\frac{2(n-1)}{n-2}}(\pa M)$};
  \node (H) [below of=T] {$H^1(M)$ };
  \node (S) [right of=H] {$L^{\frac{2n}{n-2}}(M)$};
  \draw[->] (T) to node {\small  Poisson extension} (S);
  \draw[->] (H) to node [swap] {\small  Trace} (T);
  \draw[->] (H) to node [swap] {\small Sobolev} (S);
\end{tikzpicture}}
\noindent The best constant of the Sobolev inequality plays a crucial role in solving the Yamabe problem. So does the best constant of the Sobolev trace inequality in the boundary Yamabe problem. And here, the following sharp integral inequality for the harmonic extension on the Euclidean half space
\[
\|P_{g_{\R^n}} v\|_{L^{\frac{2n}{n-2}}(\R^n_+)} \le (\Theta_{\overline B_1,g_{\R^n}})^{\frac{n-2}{2n}} \|v\|_{L^{\frac{2(n-1)}{n-2}}(\R^{n-1})}
\]
and its extremal functions (i.e., bubbles),  which is equivalent to \eqref{eq:HWY} and was proved in Hang-Wang-Yan \cite{HWY1}, play the same role in studying the variational problem \eqref{eq:theproblem}. 

However, in contrast to the other two  well studied problems, the variational problem \eqref{eq:theproblem} or \eqref{eq:energyinf} is of strong nonlocal nature. Moreover, the Euler-Lagrange equation of the functional $I[\cdot]$:
\[
v(\xi)^{\frac{n}{n-2}}=\int_M P_g(x,\xi)[(P_g v)(x)]^{\frac{n+2}{n-2}}\,\ud \mbox{vol}_g(x),
\]
which is a conformally invariant integral equation, is not the dual of any elliptic PDE. A prescribing function problem on the unit sphere, which is a Nirenberg type problem, has been studied by the second author \cite{X2}. A Kazdan-Warner type condition was obtained earlier in \cite{HWY}.

In our proof of Theorem \ref{thm:main}, the test function $v$ we use to compute  \eqref{eq:energyinf} is the cut-off of a rescaled bubble. The difficulty is to calculate the $L^{\frac{2n}{n-2}}$ norm of its Poisson extension $P_gv$. We calculate it by estimating the error between $P_gv$ and the
harmonic extension of $v$ on the Euclidean half space. Both of these two extensions are purely nonlocal. An intuitive way of estimating this error is to estimate the difference between the Poisson kernel on $M$ and 
the one on the Euclidean half space. But we are not able to accomplish it. Alternatively, we estimate the error by  
expanding it into three terms recursively, each of which satisfies a second order equation. We show that the first 
term can be calculated explicitly, and will give us a strict inequality under the assumptions of Theorem \ref{thm:main}. 
In this step, the explicit integral \eqref{eq:explicitintegral} plays a crucial role. The second term is of the same order 
as the first term. We are not able to calculate it explicitly, but we show that it is positive. If it could be calculated explicitly, then  
the dimension assumption in Theorem \ref{thm:main} might be reduced; see Remarks \ref{rem:1} and \ref{rem:2}. 
The third term is of higher order, and thus is negligible. In the low dimensional case for \eqref{eq:theproblem}, one may use more delicate test functions as in \cite{M2} and \cite{A} studying the low dimensional boundary Yamabe problem. The technical difficulty remains the same as the explicit computability of integrals involving convolutions.

This paper is organized as follows. In Section \ref{sec:idea}, we illustrate the idea of our proof.  In Section \ref{sec:nonumbilic} we consider the nonumbilic case, and in Section \ref{sec:umbilic} we study the umbilic case.
\bigskip

\noindent \textbf{Acknowledgement:} Part of this work was completed while the second named author was visiting the Department of Mathematics at the Hong Kong University of Science and Technology, to which  he is grateful   for providing  the very stimulating research environment and supports. Both authors would like to thank Professor YanYan Li for his interests and constant encouragement.

\section{Idea of the proof}\label{sec:idea}

   In the following $1\le i,j, k,l,m,p,s\le n-1$ and $1\le a ,b\le n$. For any $P\in \pa M$,  choose geodesic normal coordinates $x'=(x_1, \dots , x_{n-1})$ on the boundary centered at $P$. We say that $x=(x',x_n)$ are the Fermi coordinates of the point $\exp_{x'}(x_n\eta (x')) \in M$, where $\eta(x)$ denotes the inward unit vector normal to the boundary at $x\in M$ and $x_n\ge 0$ is small. In Fermi coordinates, we have
\[
g= g_{ij}(x)\ud x_i\ud x_j+\ud x_n^2.
\]
It was computed by Escobar \cite{E}  that
\begin{align} 
g_{ij}&=\delta_{ij} -2h_{ij}x_n-\frac{1}{3} \overline R_{ikjl}x_kx_l-2h_{ij,k}x_nx_k+(-R_{ninj}+h_{ik}h_{kj})x_n^2 +O(|x|^3),\label{eq:metric-ex1}\\
g^{ij}&=\delta_{ij} +2h_{ij}x_n+\frac{1}{3} \overline R_{ikjl}x_kx_l+2h_{ij,k}x_nx_k+(R_{ninj}+3h_{ik}h_{kj})x_n^2 +O(|x|^3),\label{eq:metric-ex2}\\
\det g&= 1-2Hx_n-\frac13 \overline R_{kl} x_kx_l-2H_{,k} x_kx_n+ (-R_{nn}-|h|^2+2H^2)x_n^2 +O(|x|^3).\label{eq:3.9}
\end{align}
Here, all coefficients are computed at $P$,   $h_{ij}$ denotes the second fundamental form with respect to the inward unit normal, $|h|^2=\sum_{ij}h_{ij}^2$, $H$ denotes the mean curvature,  and $R_{ninj}, \overline R_{ikjl}$
denote components of the full Riemannian curvature tensors of $M$ and $\pa M$, respectively. Similarly, $R_{nn}, \overline R_{kl}$
denote components of Ricci tensors of $M$ and $\pa M$, respectively.

Since $\Theta_{M,g}$ is conformally invariant, we can assume there exist conformal Fermi coordinates, i.e., Fermi coordinates $x=(x',x_n)$ centered at $P$ such that
\be\label{eq:conformalfemi}
\det g =1+O(|x|^N) \quad \mbox{in }Q_{\delta} =\{x:|x'|<\delta, 0<x_n<\delta \},
\ee
where $N$ is a positive  constant that can be chosen as large as we want, and $\delta>0$; see Marques \cite{M1}. By \eqref{eq:3.9}, as well as (3.22) in \cite{E}, we have
\be \label{eq:2}
H=H_{,k}=\overline R_{kl}=0 \quad \mbox{and} \quad  R=R_{nn}=-|h|^2 \quad \mbox{ at }0,
\ee
where $R$ is the Scalar curvature.
For any $\phi\in C^2(Q_\delta)$, we have
\[
\begin{split}
\Delta_g  \phi &=\frac{1}{\sqrt{\det g}}\pa_a (g^{ab}\sqrt{\det g} \pa_b )\phi\\&= \Delta \phi+\pa_i g^{ij}\pa _j \phi +(g^{ij}-\delta^{ij}) \pa_{ij} \phi+O(|x|^{N-1}|\nabla \phi|).
\end{split}
\]

For $0<\lda\ll\delta$, let $$v_{\lda }(x')=\left(\frac{\lda}{\lda^2+|x'|^2}\right)^{\frac{n-2}{2}} \chi_{\{|x'|<\delta\}}$$ in the above coordinates. Denote $$u_\lda=P_{g}v_\lda.$$

Let $\Omega\subset Q_{\delta}$ be a smooth domain such that $Q_{\delta/2}\subset \Omega\subset Q_{3\delta/4}$ and $\Omega$ is rotationally symmetric in the $x'$-variables. It follows that $u_\lda\ge 0$ and
\begin{align*}\begin{cases}
&L_gu_{\lda}=0 \quad \mbox{in }\Omega , \\
& u_{\lda}(x',0)=v_{\lda}(x')\quad \mbox{for }(x',0)\in \pa' \Omega,\\
& u_{\lda}=O(\lda^{\frac{n-2}{2}}) \quad \mbox{for }x\in \pa'' \Omega,
\end{cases}
\end{align*}
where $\pa' \Omega= \pa \Omega \cap \{x_n=0\}$, $\pa'' \Omega= \pa \Omega \cap \{x_n>0\}$ and the estimate of $u$ on $\pa'' \Omega$ can be proved as follows. We know there exists a conformal metric $\rho^{\frac{4}{n-2}}g$ such that $R_{\rho^{\frac{4}{n-2}}g}\equiv 0$. Then using \eqref{eq:conformalchange}, we have on $\pa'' \Omega$ that
\begin{align*}
u_\lda(x)=(P_{g}v_\lda)(x)\le C (P_{\rho^{\frac{4}{n-2}}g}v_\lda)(x)&\le C\int_{|x'|\le\delta}\frac{x_n}{(x_n^2+|x'-\xi|^2)^{\frac{n}{2}}}\frac{\lda^{\frac{n-2}{2}}}{(|\xi|^2+\lda^2)^{\frac{n-2}{2}}}\,\ud\xi\\
&\le C\lda^{\frac{n-2}{2}},
\end{align*}
where in the second inequality we used Lemma 2.2 in \cite{HWY}.

Define
\[
U_\lda(x)= \left(\frac{\lda}{|x'|^2+(x_n+\lda)^2}\right)^{\frac{n-2}{2}}
\]
which is a harmonic function in upper half space.
Let
\[
W_\lda=u_\lda-U_\lda.
\]
Then we have
\be \label{Eq:reducedequation}
\begin{cases}
&L_g W_{\lda}=F[U_\lda]\quad \mbox{in }\Omega  , \\
& W_{\lda}(x',0)=0\quad \mbox{in }\pa' \Omega,  \\
& W_{\lda}=u_{\lda}-U_{\lda} \quad \mbox{for }x\in \pa'' \Omega,
\end{cases}
\ee
where
\[
F[U_\lda]=(-L_g-\Delta)U_\lambda=\Big(\pa_i g^{ij}\pa _j +(g^{ij}-\delta^{ij}) \pa_{ij} -\frac{n-2}{4(n-1)}R_g \Big)U_\lda +O(|x|^{N-1} |\nabla U_\lda|).
\]

By direct computations, for $N>n$ we have
\begin{align}
\int_{\pa M} |v_{\lda}|^{\frac{2(n-1)}{n-2}}\,\ud s_g &= \int_{|x'|<\delta}  |v_{\lda}|^{\frac{2(n-1)}{n-2}} (1+O(|x'|^N))\,\ud x' \nonumber  \\&=2^{-(n-1)} n\w_n+O(\lda^{n-1}).
\label{eq:denominator-part}
\end{align}
We are going to show that
\begin{align}\label{eq:expansionofvolumn}
\int_{M} |u_\lda|^{\frac{2n}{n-2}}\,\ud vol_g& =  \int_{\Omega} (U_\lda+W_\lda)^{\frac{2n}{n-2}}(1+O(|x|^N))\,\ud x +O(\lda^n) \nonumber \\&
=   \int_{\Omega} \left( U_\lda ^{\frac{2n}{n-2}}  +\frac{2n}{n-2}U_\lda ^{\frac{n+2}{n-2}} W_\lda +\frac{n(n+2)}{(n-2)^2}U_\lda ^{\frac{4}{n-2}} W_\lda^2\right)\,\ud x +h.o.t.\nonumber\\&
=2^{-n} \w_n +\int_{\Omega}\left( \frac{2n}{n-2}  U_\lda ^{\frac{n+2}{n-2}} W_\lda +\frac{n(n+2)}{(n-2)^2}U_\lda ^{\frac{4}{n-2}} W_\lda^2 \right)\,\ud x +h.o.t.,
\end{align}
that is, the leading terms in the expansion of $\int_{M} |u_\lda|^{\frac{2n}{n-2}}\,\ud vol_g$ are the first three terms on the right hand side of \eqref{eq:expansionofvolumn}. The most important part of the proof will be to show
\[
\int_{\Omega}\left(   U_\lda ^{\frac{n+2}{n-2}} W_\lda +\frac{n+2}{2(n-2)}U_\lda ^{\frac{4}{n-2}} W_\lda^2 \right)\,\ud x >0.
\]
To verify the above inequality, the next crucial step is to solve the leading terms of $W_{\lda}$, which we will divide it into two cases in the following two sections. The explicit integral \eqref{eq:explicitintegral} plays a crucial role in our calculations.

\section{Nonumbilic boundary case}\label{sec:nonumbilic}

Suppose the $P\in \pa M$ is a nonumbilic point. There exist conformal Fermi coordinates centered at $P$ as in previous section.    By \eqref{eq:metric-ex1}, \eqref{eq:metric-ex2} and \eqref{eq:2}, we have
\begin{align*}
&F[U_\lda](x)\\
&=(\pa_i g^{ij}\pa _j +(g^{ij}-\delta^{ij}) \pa_{ij} -\frac{n-2}{4(n-1)}R_g )U_\lda+O(|x|^{N-1} |\nabla U_\lda|)\\&
=2h_{ij,i} x_n\pa_j U_{\lda}+[2h_{ij}x_n+\frac{1}{3} \overline R_{ikjl}x_kx_l+2h_{ij,k}x_nx_k+(R_{ninj}+3h_{ik}h_{kj})x_n^2]\pa_{ij}U_{\lda}\\ &
\quad +\frac{n-2}{4(n-1)} |h|^2 U_{\lda}+ O(|x| U_{\lda}),
\end{align*}
where
\[
\pa_j U_\lda= -(n-2)\lda^{\frac{n-2}{2}}\Big((\lda+x_n)^2+|x'|^2\Big)^{-\frac{n}{2}}x_j
\]
and
\[
\pa_{ij} U_{\lda}=  -(n-2)\lda^{\frac{n-2}{2}}\Big((\lda+x_n)^2+|x'|^2\Big)^{-\frac{n}{2}}\delta_{ij}+n(n-2)\lda^{\frac{n-2}{2}}\Big((\lda+x_n)^2+|x'|^2\Big)^{-\frac{n+2}{2}}x_ix_j.
\]

Denote the solution of the Dirichlet problem
\[
-\Delta u= f \quad \mbox{in }\Omega, \quad u=0 \quad \mbox{on }\pa \Omega
\]
by
\be \label{eq:green-represent}
u= \mathcal{G}(f)=\int_{\Omega} G(x,y) f(y)\,\ud y,
\ee where $G(x,y)$ is the Green's function. We know that $G(x,y)=G(y,x)$ and there exists $C>0$ such that $G(x,y)\le C|x-y|^{2-n}$.

Define
\[
W^{(1)}_{\lda}= \mathcal{G}(F[U_\lda]) \quad  \mbox{and}\quad  W^{(2)}_{\lda}= \mathcal{G}(F(W^{(1)}_{\lda})).
\]
We assume $n\ge 10$ in the following because of our assumption on the dimension in our main theorem.

The next three lemmas show that $W^{(1)}_{\lda}$ and $W^{(2)}_{\lda}$ are the leading terms in $W_{\lda}$, and the rest error will be a higher order term.
\begin{lem}\label{lem:leading}
There exists $C>0$ such that for all $s=1,2,3, 4$, and all $x\in\Omega$,
\begin{equation}\label{eq:gradientestimate}
|W_\lda^{(1)}(x)|+|x+\lda e_n|^s |\nabla^s W_\lda^{(1)}(x)|\le C \lda^{\frac{n-2}{2}} |x+\lda e_n|^{3-n}.
\end{equation}
\end{lem}
\begin{proof}
For $s=1,2,3, 4$, it is direct to see that
\begin{equation}\label{eq:errorFU}
|F[U_\lda](x)| +|x+\lda e_n|^s|\nabla^s (F[U_\lda](x))| \le  C \lda^{\frac{n-2}{2}} |x+\lda e_n|^{-n+1} \quad\mbox{ in } \Omega.
\end{equation}
for some constant $C$ independent of $\lda$. Then  we have for all $x\in \Omega$
\begin{equation}\label{eq:solveleading}
\begin{split}
|W_\lda^{(1)}(x)|\le C \int_{\Omega}\frac{|F[U_\lambda](y)|}{|x-y|^{n-2}}\,\ud y&\le C \lda^{\frac{n-2}{2}} \int_{\Omega}\frac{1}{|x-y|^{n-2}|y+\lda e_n|^{n-1}}\,\ud y\\
&\le
C \lda^{\frac{n-2}{2}} |x+\lda e_n|^{3-n}\quad\mbox{if }n> 3.
\end{split}
\end{equation}

Now we estimate $|\nabla^s W_\lda^{(1)}(x)|$. By standard elliptic estimates, we have that
\[
|\nabla^s W_\lda^{(1)}(x)|\le C \lda^{\frac{n-2}{2}} \quad\mbox{ in } \Omega\setminus Q_{\delta/8}.
\]
To obtain the estimate in $Q_{\delta/8}$, we use a scaling argument. Let $x\in Q_{\delta/8}$ be arbitrarily fixed.  Let $r=|x+\lda e_n|/2$, and
\[
v(y)=r^{n-3}W_\lda^{(1)}(x+ry).
\]
where
\[
y\in \mathcal{O}:=\{y\in B_1: x+ry\in \Omega\}.
\]
Then $\overline B^+_1\subset\mathcal{O}$ (because $x\in Q_{\delta/8}$), and
\[
-\Delta v=r^{n-1}F[U_\lda](x+ry)=:f(y)\quad\mbox{in }\mathcal{O}.
\]
Also, it follows from \eqref{eq:solveleading} and \eqref{eq:errorFU} that
\[
|v|+|f|+|\nabla^s f|\le C\lda^{\frac{n-2}{2}}\quad\mbox{in }\mathcal{O}.
\]
Case 1: if $B_{1/8}\subset\mathcal{O}$, then it follows from the interior estimates for the Poisson equations that
\[
|\nabla^s v(0)|\le C\lda^{\frac{n-2}{2}}.
\]
Case 2: if $B_{1/8}\not\subset\mathcal{O}$, since $v=0$ on $\mathcal{O}\cap \partial\R^n_+$, it follows from the global estimates up to the boundary for the Poisson equations that
\[
|\nabla^s v(0)|\le C\lda^{\frac{n-2}{2}}.
\]
Rescaling to $W_\lda^{(1)}$, we obtain \eqref{eq:gradientestimate}.
\end{proof}
\begin{lem}\label{lem:leading2}
There exists $C>0$ such that for both $s=1,2$, and all $x\in\Omega$,
\begin{equation}\label{eq:gradientestimate2}
|W_\lda^{(2)}(x)|+|x+\lda e_n|^s |\nabla^s W_\lda^{(2)}(x)|\le C \lda^{\frac{n-2}{2}} |x+\lda e_n|^{4-n}.
\end{equation}
\end{lem}
\begin{proof}
It follows from Lemma \ref{lem:leading} that
\[
|F[W^{(1)}_\lambda(x)]|+|x+\lda e_n|^s|\nabla^s (F[W^{(1)}_\lambda](x))| \le C \lda^{\frac{n-2}{2}} |x+\lda e_n|^{2-n},
\]
which implies as before that
\begin{equation}\label{eq:solvesecond}
\begin{split}
|W_\lda^{(2)}(x)|\le C \int_{\Omega}\frac{|F[W^{(1)}_\lambda](y)|}{|x-y|^{n-2}}\,\ud y&\le C \lda^{\frac{n-2}{2}} \int_{\Omega}\frac{1}{|x-y|^{n-2}|y+\lda e_n|^{n-2}}\,\ud y\\[2mm]
&\le
C \lda^{\frac{n-2}{2}} |x+\lda e_n|^{4-n}\quad\mbox{if }n> 4.
\end{split}
\end{equation}
The rest of the proof is the same as that of Lemma \ref{lem:leading}.
\end{proof}

\begin{lem}\label{lem:leading3}
Let $W^{(3)}_{\lda}:=W_{\lda}-W^{(1)}_{\lda}-W^{(2)}_{\lda} $. There exists $C>0$ such that for all $x\in\Omega$,
\begin{equation}\label{eq:solveleast}
|W_\lda^{(3)}(x)|\le C \lda^{\frac{n-2}{2}} |x+\lda e_n|^{5-n}.
\end{equation}
\end{lem}
\begin{proof}
 By \eqref{Eq:reducedequation}, as well as the definitions of $W_{\lda}^{(1)}$ and $W_{\lda}^{(2)}$, we know that $W_{\lda}^{(3)}$ satisfies
\[
\begin{split}
&L_g W^{(3)}_{\lda}=F[W^{(2)}_{\lda}] \quad \mbox{in }\Omega  , \\
& W^{(3)}_{\lda} (x',0)=0,\\
& W^{(3)}_{\lda}=u_{\lda}-U_{\lda} \quad \mbox{for }x\in \pa'' \Omega.
\end{split}
\]
We decompose $W^{(3)}_{\lda}=W^{(31)}_{\lda}+W^{(32)}_{\lda}$, where
\[
\begin{split}
&L_g W^{(31)}_{\lda}=F[W^{(2)}_{\lda}] \quad \mbox{in }\Omega  , \\
& W^{(31)}_{\lda}=0 \quad \mbox{for }x\in \pa \Omega,
\end{split}
\]
and
\[
\begin{split}
&L_g W^{(32)}_{\lda}=0\quad \mbox{in }\Omega  , \\
& W^{(32)}_{\lda} (x',0)=0,\\
& W^{(32)}_{\lda}=u_{\lda}-U_{\lda} \quad \mbox{for }x\in \pa'' \Omega.
\end{split}
\]
Denote the first Dirichlet eigenvalue of $L_{g}$ in $\om$ by $\lambda_1 (L_g,\Omega) $. Since $\om$ is a submanifold of $M$,  $\lambda_1 (L_g,\Omega) \ge \lambda_1(L_{ g})>0$, the Green function for the equation of $W^{(31)}_\lambda$ exists, and  we have for all $x\in \Omega$
\begin{equation}\label{eq:solveleast1}
\begin{split}
|W_\lda^{(31)}(x)|\le C \int_{\Omega}\frac{|F[W^{(2)}_{\lda}](y)|}{|x-y|^{n-2}}\,\ud y&\le C \lda^{\frac{n-2}{2}} \int_{\Omega}\frac{1}{|x-y|^{n-2}|y+\lda e_n|^{n-3}}\,\ud y\\
&\le
C \lda^{\frac{n-2}{2}} |x+\lda e_n|^{5-n}\quad\mbox{if }n> 5.
\end{split}
\end{equation}
Again, since $\lambda_1(L_g,\Omega)>0$, we have by the comparison principle,
\begin{equation}\label{eq:solveleast2}
|W^{(32)}_{\lda}(x)|\le C\sup_{x\in \pa'' \Omega}|u_{\lda}-U_{\lda}|\le C \lda^{\frac{n-2}{2}}.
\end{equation}
Indeed, let $\rho$ be the solution of
\[
L_g\rho=0\quad\mbox{in }\Omega,\quad\rho=1\quad\mbox{on }\partial\Omega.
\]
Since $\lambda_1(L_g,\Omega)>0$, we have that $\rho>0$ in $\overline\Omega$. $\forall~\varepsilon>0$, let $\varphi=(\sup_{x\in \pa'' \Omega}|u_{\lda}-U_{\lda}|+\varepsilon)\rho\pm W^{(32)}_{\lda}$. Then
\[
L_g\varphi=0\quad\mbox{in }\Omega,\quad\varphi\ge\varepsilon\quad\mbox{on }\partial\Omega.
\]
Thus, $\varphi>0$ in $\overline\Omega$. By sending $\varepsilon\to 0$, we obtain \eqref{eq:solveleast2}.

The conclusion follows from \eqref{eq:solveleast1} and \eqref{eq:solveleast2}.
\end{proof}

\begin{proof}[Proof of Theorem \ref{thm:main} (i)]
As explained in \eqref{eq:expansionofvolumn}, we are going to calculate
\[
\int_{\Omega}U_\lda ^{\frac{n+2}{n-2}} W_\lda \,\ud x\quad\mbox{and}\quad \int_{\Omega}U_\lda ^{\frac{4}{n-2}} W_\lda^2 \,\ud x.
\]

It follows from Lemmas \ref{lem:leading}, \ref{lem:leading2} and \ref{lem:leading3} that
\begin{align}
\int_{\Omega} U_{\lda}^{\frac{n+2}{n-2}} W_{\lda}&= \int_{\Omega} U_{\lda}^{\frac{n+2}{n-2}} W^{(1)}_{\lda}+ \int_{\Omega} U_{\lda}^{\frac{n+2}{n-2}} W^{(2)}_{\lda}+ \int_{\Omega} U_{\lda}^{\frac{n+2}{n-2}} W^{(3)}_{\lda}\nonumber\\
&= \int_{\Omega} U_{\lda}^{\frac{n+2}{n-2}} W^{(1)}_{\lda}+ \int_{\Omega} U_{\lda}^{\frac{n+2}{n-2}} W^{(2)}_{\lda}+ O(\lda^{3})\nonumber\\
&= \int_{\Omega} \mathcal{G}( U_{\lda}^{\frac{n+2}{n-2}}) F[U_{\lda}]+\int_{\Omega}  \mathcal{G}( U_{\lda}^{\frac{n+2}{n-2}}) F[W^{(1)}_{\lda}]+O(\lda^{3}),\label{eq:decomp1}
\end{align}
where we used Fubini theorem and $G(x,y)=G(y,x)$ in the first two terms of the last equality.

\emph{Step 1:} We first calculate
\[
\int_{\Omega} \mathcal{G}( U_{\lda}^{\frac{n+2}{n-2}}) F[U_{\lda}].
\]

Since $U_{\lda}$ is radial symmetric in $x'$ and so is $\Omega$, $\mathcal{G}( U_{\lda}^{\frac{n+2}{n-2}})$  depends only on $|x'| $ and $x_n $.
By the symmetry in $x'$, it follows from \eqref{eq:2} that
\begin{align*}
&\int_{\Omega} \mathcal{G}( U_{\lda}^{\frac{n+2}{n-2}})F[U_{\lda}] \\&
= \int_{\Omega}\mathcal{G}( U_{\lda}^{\frac{n+2}{n-2}})\Big( ( R_{ninj}+3h_{ik}h_{kj})x_n^2\pa_{ij}U_{\lda}+\frac{n-2}{4(n-1)} |h|^2   U_{\lda} +O(|x|U_\lda)\Big)\\
&= \lda^{\frac{n-2}{2}} |h|^2  \int_{\Omega} \mathcal{G}( U_{\lda}^{\frac{n+2}{n-2}}) f_\lda(x) + O\left(\int_{\Omega}\mathcal{G}( U_{\lda}^{\frac{n+2}{n-2}})|x|U_\lda\right)\\
&=\lda^{2} |h|^2  \int_{\Omega_\lda} \mathcal{G}_\lda( U_{1}^{\frac{n+2}{n-2}}) f_1(x) + O\left(\lda^3\int_{\Omega_\lda}\mathcal{G}_\lda( U_{1}^{\frac{n+2}{n-2}})|x|U_1\right),
\end{align*}
where $\Omega_\lambda=\{x/\lambda: x\in\Omega\}$,
\[
\mathcal{G}_{\lda}(f)=\int_{\Omega_\lambda} G_\lambda(x,y) f(y)\,\ud y \quad \mbox{with}\quad G_\lambda(x,y)=\lambda^{n-2}G(\lambda x, \lambda y),
\]
and
\begin{align*}
f_\lda(x):&= -2(n-2) x_n^2((\lda+x_n)^2+|x'|^2)^{-\frac{n}{2}}\\&
\quad + \frac{2n(n-2)}{n-1}x_n^2 |x'|^2 ((\lda+x_n)^2+|x'|^2)^{-\frac{n+2}{2}}+\frac{n-2}{4(n-1)}((\lda+x_n)^2+|x'|^2)^{-\frac{n-2}{2}}.
\end{align*}

We are going to estimate $\mathcal{G}_\lda( U_{1}^{\frac{n+2}{n-2}})$. Let
\[
\overline G(x,y)=\frac{1}{n(n-2)\w_n}\left(\frac{1}{|x-y|^{n-2}}-\frac{1}{|(x',-x_n)-y|^{n-2}}\right) \quad \mbox{for } x_n>0
\]
be the Green's function of the Poisson equation on the upper half space. We have
\begin{align}\label{eq:explicitintegral}
\int_{\R^n_+} \overline G(x,y) U_{1}(y)^{\frac{n+2}{n-2}}\,\ud y
=\frac{1}{2n}   x_n((1+x_n)^2+|x'|^2)^{-\frac{n}{2}}
\end{align}
by observing that
\[
-\Delta\left(\frac{1}{2n}   x_n((1+x_n)^2+|x'|^2)^{-\frac{n}{2}}\right)=U_{1}(y)^{\frac{n+2}{n-2}}.
\]
In particular, $\int_{\R^n_+} \overline G(x,y) U_{1}(y)^{\frac{n+2}{n-2}}\,\ud y  \le C x_n \lda^{n}$ on $\pa \Omega_{\lda}$. Since
\[
\Delta \left( \int_{\R^n_+} \overline G(x,y) U_{1}(y)^{\frac{n+2}{n-2}}\,\ud y - \int_{\Omega_{\lda}}G_{\lda} (x,y)U_{1}(y)^{\frac{n+2}{n-2}}\,\ud y  \right)=0 \quad \forall ~x\in \Omega_{\lda},
\]
it follows from the maximum principle that in $\Omega_\lda$,
\[
0\le \int_{\R^n_+} \overline G(x,y) U_{1}(y)^{\frac{n+2}{n-2}}\,\ud y - \int_{\Omega_{\lda}}G_{\lda} (x,y)U_{1}(y)^{\frac{n+2}{n-2}}\,\ud y \le C\lda^{n}  x_n.
\]
Hence, we have that in $\Omega_\lda$,
\begin{equation}\label{eq:theleadingterm}
\mathcal{G}_\lda (U_{1}^{\frac{n+2}{n-2}})=\frac{1}{2n} x_n |x+e_n|^{-n}+\lda^n x_nO(1).
\end{equation}
Consequently, as long as $n\ge5$, there holds
\[
\int_{\Omega_\lda}\mathcal{G}_\lda( U_{1}^{\frac{n+2}{n-2}})|x|U_1\le C.
\]

Therefore, we have
\begin{align*}
\int_{\Omega} \mathcal{G}( U_{\lda}^{\frac{n+2}{n-2}})F[U_{\lda}] 
&=\frac{\lda^{2} |h|^2}{2n}\int_{\Omega_{\lda}}  y_n((1+y_n)^2+|y'|^2)^{-\frac{n}{2}} f_1(y)\,\ud y+O(\lda^{n-1})+O(\lda^3)\\&
=\frac{\lda^{2} |h|^2}{2n}\left(-2(n-2) I_1+  \frac{2n(n-2)}{n-1} I_2+\frac{n-2}{4(n-1)} I_3\right) + O(\lda^{3}),
\end{align*}
where
\begin{align*}
I_1&= \int_{\R^n_+} y_n^3 |y+e_n|^{-2n}\,\ud y,\\
I_2&=  \int_{\R^n_+} y_n^3|y'|^2 |y+e_n|^{-2n-2}\,\ud y,\\
I_3&=   \int_{\R^n_+} y_n |y+e_n|^{-2n+2}\,\ud y.
\end{align*}
If $n> 4$, by changing variables $y'=(1+y_n)z'$ we have

\begin{align*}
I_1&= \int_{\R^{n}_+} y_n^3 |y+e_n|^{-2n}\,\ud y\\
&= \int_0 ^\infty  \frac{y_n^3}{(1+y_n)^{n+1}}\,\ud y_n \int_{\R^{n-1}} (1+|z'|^2)^{-n}\,\ud z' \\&
=\frac{6|\mathbb S^{n-2}| }{n(n-1)(n-2)(n-3)} \int_0^\infty \frac{r^{n-2}}{(1+r^2)^{n}}\,\ud r,
\end{align*}
\begin{align*}
I_2&= \int_{\R^{n}_+} y_n^3 |y'|^2 |y+e_n|^{-2n-2}\,\ud y +O(\lda)\\
&= \int_0 ^\infty  \frac{y_n^3}{(1+y_n)^{n+1}}\,\ud y_n \int_{\R^{n-1}} |z'|^2(1+|z'|^2)^{-n-1}\,\ud z' \\&
=\frac{6|\mathbb S^{n-2}| }{n(n-1)(n-2)(n-3)} \int_0^\infty \frac{r^{n}}{(1+r^2)^{n+1}}\,\ud r\\&
=\frac{6|\mathbb S^{n-2}|}{n(n-1)(n-2)(n-3)} \frac{n-1}{2n}\int_0^\infty \frac{r^{n-2}}{(1+r^2)^{n}}\,\ud r,
\end{align*}
and
\begin{align*}
I_3&= \int_{\R^{n}_+} y_n |y+e_n|^{-2n+2}\,\ud y\\
&= \int_0 ^\infty  \frac{y_n}{(1+y_n)^{n-1}}\,\ud y_n \int_{\R^{n-1}} (1+|z'|^2)^{-n+1}\,\ud z'\\&
=\frac{|\mathbb S^{n-2}|}{(n-2)(n-3)} \int_0^\infty \frac{r^{n-2}}{(1+r^2)^{n-1}}\,\ud r\\
&=\frac{2|\mathbb S^{n-2}|}{(n-2)(n-3)} \int_0^\infty \frac{r^{n-2}}{(1+r^2)^{n}}\,\ud r.
\end{align*}
It follows that
\begin{align*}
-2(n-2) I_1+  \frac{2n(n-2)}{n-1} I_2+\frac{n-2}{4(n-1)} I_3 =\frac{|\mathbb S^{n-2}| (n-12)}{2n(n-1)(n-3)}\int_0^\infty \frac{r^{n-2}}{(1+r^2)^{n}}\,\ud r,
\end{align*}
which is positive if $n>12$.  Therefore,
\begin{align}\label{eq:estimateofleadingterm1}
\int_{\Omega} \mathcal{G}( U_{\lda}^{\frac{n+2}{n-2}}) F[U_{\lda}]=  \lda^2 |h|^2 \frac{|\mathbb S^{n-2}| (n-12)}{4n^2(n-1)(n-3)}\int_0^\infty \frac{r^{n-2}}{(1+r^2)^{n}}\,\ud r +O(\lda^3).
\end{align}


\emph{Step 2:} Next we calculate
\[
\int_{\Omega}  \mathcal{G}( U_{\lda}^{\frac{n+2}{n-2}}) F[W^{(1)}_{\lda}].
\]

It follows from \eqref{eq:gradientestimate} that
\begin{align*}
\int_{\Omega}  \mathcal{G}( U_{\lda}^{\frac{n+2}{n-2}}) F[W^{(1)}_{\lda}]&=\int_{\Omega}  \mathcal{G}( U_{\lda}^{\frac{n+2}{n-2}}) \Big(2h_{ij}x_n\partial_{ij}W^{(1)}_{\lda}+O(\lda^{\frac{n-2}{2}}|x+\lda e_n|^{3-n})\Big)\\
& =2\int_{\Omega}  \mathcal{G}( U_{\lda}^{\frac{n+2}{n-2}}) h_{ij}x_n\partial_{ij}W^{(1)}_{\lda}+O(\lambda^3).
\end{align*}
It follows from \eqref{eq:theleadingterm} that
\begin{equation}\label{eq:theleadingtermscaled}
\mathcal{G} (U_{\lda}^{\frac{n+2}{n-2}})(x)=\lda^{\frac{2-n}{2}}\mathcal{G}_\lda (U_{1}^{\frac{n+2}{n-2}})(\lda^{-1}x)=\frac{\lda^{n/2}}{2n} x_n |x+\lda e_n|^{-n}+\lda^{n/2}x_n O(1).
\end{equation}
Hence, using integration by parts, $H=\sum_{i} h_{ii}=0$ and $F(U_\lda)=\mathcal{G}(h_{kl}y_n\partial_{kl}U_\lda+O(U_\lda))$, we have
\begin{align*}
&2\int_{\Omega}  \mathcal{G}( U_{\lda}^{\frac{n+2}{n-2}}) h_{ij}x_n\partial_{ij}W^{(1)}_{\lda}\\&=\frac{1}{n}\lda^{n/2}\int_{\Omega} x_n^2|x+\lda e_n|^{-n}  h_{ij}\partial_{ij}W^{(1)}_{\lda}+O(\lda^{n-1})\\
&=\lda^{n/2}(n+2)\int_{\Omega}x_n^2|x+\lda e_n|^{-n-4}h_{ij}x_ix_j\mathcal{G}(F(U_\lambda))+O(\lda^{n-1})\\
&=2\lda^{n/2}(n+2)\int_{\Omega}x_n^2|x+\lda e_n|^{-n-4}h_{ij}x_ix_j\mathcal{G}(h_{kl}y_n\partial_{kl}U_\lda)+O(\lda^3)\\
&=2\lda^{2}(n+2)\int_{\Omega_\lda}x_n^2|x+ e_n|^{-n-4}h_{ij}x_ix_j\mathcal{G}_\lda(h_{kl}y_n\partial_{kl}U_1)+O(\lda^3),
\end{align*}
where we changed variables in the last identity, and the integrals on the boundary $\partial\Omega$ coming from the integration by part are absorbed in $O(\lda^{n-1})$. 

We are going to estimate $\mathcal{G}_\lda(x_n\partial_{kl}U_1)$. Denote
\[
\overline{\mathcal{G}} (f)=\int_{\R^n_+}\overline G(x,y)f(y)\,\ud y.
\]
Then
\begin{equation}\label{eq:auxupperbound}
|\overline{\mathcal{G}}( y_n\pa_{ij} U_{1} )(x)|\le C\int_{\R^n_+}\frac{1}{|x-y|^{n-2}|y+ e_n|^{n-1}}\,\ud y\le C|x+e_n|^{3-n}.
\end{equation}
By the maximum principle for harmonic functions, we have
\begin{equation}\label{eq:aux211}
\|\mathcal{G}_\lambda (y_n\pa_{ij} U_{1} )-\overline{\mathcal{G}}(y_n\pa_{ij} U_{1} )\|_{L^\infty(\Omega_\lambda)}\le \|\overline{\mathcal{G}}( x_n\pa_{ij} U_{1} )\|_{L^\infty(\partial\Omega_\lambda)}\le C\lambda^{n-3}.
\end{equation}
Therefore,
\begin{align*}
&\int_{\Omega}  \mathcal{G}( U_{\lda}^{\frac{n+2}{n-2}}) F[W^{(1)}_{\lda}]\\
&=2\lda^{2}(n+2)\int_{\Omega_\lda}x_n^2|x+ e_n|^{-n-4}h_{ij}x_ix_j\overline{\mathcal{G}}(h_{kl}y_n\partial_{kl}U_1)+O(\lda^3)\\
&=2\lda^{2}(n+2)\int_{\R^n_+}x_n^2|x+ e_n|^{-n-4}h_{ij}x_ix_j\overline{\mathcal{G}}(h_{kl}y_n\partial_{kl}U_1)+O(\lda^3)\\
&= 2\lda^{2}n(n^2-4)\int_{\R^n_+}x_n^2|x+ e_n|^{-n-4}h_{ij}x_ix_j\overline{\mathcal{G}}(h_{kl}y_ky_ly_n|y+ e_n|^{-2-n})+O(\lda^3).
\end{align*}

 Now we calculate  $\Phi(x):=\overline{\mathcal{G}} ( h_{ij}y_iy_j y_n |y+e_n|^{-2-n})(x)$. It solves
\[
-\Delta \Phi= h_{ij} x_ix_j \eta(|x'|,x_n) \quad \mbox{in }\R^n_+, \quad \Phi= 0  \quad \mbox{on }\pa \R^n_+,
\]
where $\eta(|x'|,x_n)= x_n (|x'|^2+ (x_n+1)^2)^{-\frac{2+n}{2}}$.   We look for the unique function $\Phi$ in a form of $$\Phi= h_{ij}x_ix_j V(|x'|,x_n)$$ for some $V=V(r,s)$. Then $V$ satisfies
\begin{equation}\label{eq:highlap-0}
-\partial_{rr} V-\frac{n+2}{r}\partial_rV-\partial_{ss}V= \eta(r,s), \quad V(r,0)= 0\quad \mbox{for }r,s>0,
\end{equation}
where $\sum_i h_{ii}=0$ is used. The operator on the left hand side of \eqref{eq:highlap-0} can be considered as the Laplacian operator $-\Delta$ in $\R^{n+4}$ applying to $V$ which is radial in the first $n+3$ coordinates. Therefore, for $\tilde V(z)=V(|z'|,z_{n+4})$ with $z'=(z_1,\cdots,z_{n+3})$, it satisfies
\begin{equation}\label{eq:highlap-1}
-\Delta \tilde V(z)= z_{n+4} |z+ e_{n+4}|^{-2-n}  \quad \mbox{in }\R^{n+4}_+, \quad \tilde V= 0  \quad \mbox{on }\pa\R^{n+4}_+.
\end{equation}
This $\tilde V$ can be solved using Green's function of $-\Delta$ in $\R^{n+4}_+$, and thus, $\tilde V>0$ everywhere.  Moreover, by the same argument as of \eqref{eq:solveleading}, we have that
\[
\tilde V(z)\le C   |z+ e_{n+4}|^{1-n} \quad \mbox{in }\R^{n+4}_+.
\]
Hence, we have for all $r,s\ge 0$,
\[
0<V(r,s)\le C  (r^2+(1+s)^2)^\frac{1-n}{2}.
\]

Therefore,
\[
\int_{\Omega}  \mathcal{G}( U_{\lda}^{\frac{n+2}{n-2}}) F[W^{(1)}_{\lda}]=2n(n^2-4)\lambda^{2}\int_{\R^n_+}x_n^2  |x+e_n|^{-(n+4)} V(|x'|,x_n)(h_{ij}x_ix_j)^2 +O(\lambda^3).
\]
Since $\sum_{i} h_{ii}=0$, we have
\begin{align*}
\int_{\mathbb S^{n-2}}(h_{ij}x_ix_j)^2 &= \int_{\mathbb S^{n-2}} \sum_{i}h_{ii}^2x_i^4+2\sum_{i\neq j}h_{ij}^2x_i^2x_j^2+\sum_{i\neq j}h_{ii}h_{jj}x_i^2x_j^2\\
&=\int_{\mathbb S^{n-2}} x_1^4\Big(\sum_{i}h_{ii}^2+\frac{2}{3}\sum_{i\neq j}h_{ij}^2+\frac 13\sum_{i\neq j}h_{ii}h_{jj}\Big)\\
&=\frac{2}{3} |h|^2\int_{\mathbb S^{n-2}} x_1^4,
\end{align*}
where in the second equality we used
\begin{equation}\label{eq:trick}
\int_{\mathbb S^{n-2}}x_1^2x_2^2=\frac 13 \int_{\mathbb S^{n-2}}x_1^4=\frac{|\mathbb S^{n-2}|}{(n-1)(n+1)}.
\end{equation}

Thus,
\begin{align}
\int_{\Omega}  \mathcal{G}( U_{\lda}^{\frac{n+2}{n-2}}) F[W^{(1)}_{\lda}]&=\frac{4n(n^2-4)}{3}\lambda^{2}|h|^2\int_{\R^n_+}x_n^2  |x+e_n|^{-(n+4)} V(|x'|,x_n)x_1^4+O(\lambda^3)\nonumber\\
&=C(n)|h^2|\lda^2+ O(\lambda^3)\label{eq:estimateofleadingterm2}
\end{align}
for some $C(n)>0$.

\emph{Step 3:} Finally, we estimate
\[
\int_{\Omega}U_\lda ^{\frac{4}{n-2}} W_\lda^2 \,\ud x.
\]

It follows from \eqref{eq:solveleading}, \eqref{eq:solvesecond} and \eqref{eq:solveleast} that
\[
\int_{\Omega} U_\lda ^{\frac{4}{n-2}} W_\lda^2 = \int_{\Omega} U_\lda ^{\frac{4}{n-2}} (W^{(1)}_\lda)^2 +O(\lda^3).
\] By the definition of $W^{(1)}_\lda$, for $n> 5$ we have
\begin{align*}
W^{(1)}_\lda& = 2\mathcal{G} (h_{ij}x_n\pa_{ij} U_{\lda} +O(U_\lda))\\&
= 2\mathcal{G} ( h_{ij}  x_n\pa_{ij} U_{\lda} )+O(\lda^{\frac{n-2}{2}}|x+\lda e_n|^{4-n}).
\end{align*}
Therefore, using \eqref{eq:auxupperbound} and \eqref{eq:aux211}, we have
\begin{align}
\int_{\Omega} U_\lda ^{\frac{4}{n-2}} W_\lda^2&= 4\lambda^2 \int_{\Omega_\lambda} U_1 ^{\frac{4}{n-2}} (\mathcal{G}_\lambda (h_{ij}   x_n\pa_{ij} U_{1} ))^2 + O(\lda^3)\nonumber\\
&= 4\lambda^2 \int_{\Omega_\lambda} U_1 ^{\frac{4}{n-2}} (\overline{\mathcal{G}} ( h_{ij}  x_n \pa_{ij} U_{1} ))^2 + O(\lda^3)\nonumber\\
&= 4\lambda^2 \int_{\R^n_+} U_1 ^{\frac{4}{n-2}} \big(\overline{\mathcal{G}} ( x_n h_{ij}  \pa_{ij} U_{1} )\big)^2 + O(\lda^3)\nonumber\\
&=4\lambda^2n^2(n-2)^2 \int_{\R^n_+} U_1 ^{\frac{4}{n-2}} \big(\overline{\mathcal{G}} ( x_n h_{ij}x_ix_j  |x+e_n|^{-2-n})\big)^2 + O(\lda^3)\nonumber\\
&=4n^2(n-2)^2\lambda^2 \int_{\R^n_+} U_1 ^{\frac{4}{n-2}} V^2(|x'|,x_n) (h_{ij}  x_ix_j)^2 + O(\lda^3)\nonumber \\
&=\frac{8n^2(n-2)^2}{3}\lambda^2 |h|^2\int_{\R^n_+} U_1 ^{\frac{4}{n-2}} V^2(|x'|,x_n) x_1^4+ O(\lda^3)\nonumber\\
&=C(n)|h|^2\lda^2+ O(\lda^3) \label{eq:estimateofleadingterm3}
\end{align}
for some $C(n)>0$, where $V$ is the same one as in \eqref{eq:highlap-0}.

Consequently, it follows from \eqref{eq:decomp1}, \eqref{eq:estimateofleadingterm1}, \eqref{eq:estimateofleadingterm2} and \eqref{eq:estimateofleadingterm3} that
\begin{align}
&\int_{M} |u_\lda|^{\frac{2n}{n-2}}\,\ud vol_g \nonumber \\& =  \int_{\Omega} (U_\lda+W_\lda)^{\frac{2n}{n-2}}(1+O(|x|^N))\,\ud x +O(\lda^3) \nonumber \\&
=   \int_{\Omega} \left( U_\lda ^{\frac{2n}{n-2}}  +\frac{2n}{n-2}U_\lda ^{\frac{n+2}{n-2}} W_\lda +\frac{n(n+2)}{(n-2)^2}U_\lda ^{\frac{4}{n-2}} W_\lda^2 +O(|W_\lda|^{\frac{2n}{n-2}})\right)\,\ud x  +O(\lda^3)\nonumber \\&
\ge 2^{-n} \w_n +\frac{2n}{n-2} \int_{\Omega}U_\lda ^{\frac{n+2}{n-2}} W_\lda\,\ud x+\frac{n(n+2)}{(n-2)^2}U_\lda ^{\frac{4}{n-2}} W_\lda^2 +O(\lda^{\frac{2n}{n-2}}) \nonumber \\&
\ge 2^{-n} \w_n +C(n)\lda^2 |h|^2 +O(\lda^{\frac{2n}{n-2}})\label{eq:numerate1}
\end{align}
for some $C(n)>0$ as long as $n\ge 12$.

By \eqref{eq:denominator-part} and \eqref{eq:numerate1}, we have
\begin{align*}
\frac{\int_{M} |u_\lda|^{\frac{2n}{n-2}}\,\ud vol_g}{(\int_{\pa M} |v_{\lda}|^{\frac{2(n-1)}{n-2}}\,\ud s_g )^{\frac{n}{n-1}}}\ge \frac{2^{-n} \w_n +C(n)\lda^2|h|^2 +O(\lda^{\frac{2n}{n-2}})}{(2^{-(n-1)} n\w_n+O(\lda^{n-1}))^{\frac{n}{n-1}}}>  n^{-\frac{n}{n-1}} \w_n^{-\frac{1}{n-1}},
\end{align*}
if $n\ge 12$.
\end{proof}






\begin{rem} \label{rem:1}
From Steps 2 and 3 in the above, if one can explicitly calculate $V$ in \eqref{eq:highlap-0}, or equivalently $\tilde V$ in \eqref{eq:highlap-1}, then the dimension assumption might be reduced.

\end{rem}

\section{Umbilic boundary case}\label{sec:umbilic}

Suppose $ \pa M$ is umbilic. As in the previous section, we assume that there exists a conformal Fermi coordinates $x=(x',x_n)$ centered at some point  $P\in \pa M$. Below we collect  some facts which were proved by Marques \cite{M1} (see also Lemma 2.3 in \cite{A}).
First,
\be
h_{ij}=O(|x|^{N-1}),
\ee
and
\be
\begin{split}
g^{ij}=& \delta_{ij}+\frac13 \overline R_{ikjl} x_kx_l+R_{ninj} x_n^2+\frac16 \overline R_{ikjl,m}x_{k}x_lx_m+R_{ninj,k}x_n^2 x_k +\frac13 R_{ninj,n} x_n^3\\&
+(\frac{1}{20} \overline R_{ikjl,mp}+\frac{1}{15} \overline R_{iksl} \overline R_{jmsp})x_{k}x_l x_mx_p\\
&+(\frac12 R_{ninj,kl} +\frac13 Sym_{ij}(\overline R_{iksl} R_{nsnj}))x_n^2 x_k x_l\\&
+\frac13 R_{ninj,nk} x_n^3x_k + (\frac{1}{12} R_{ninj,nn}+\frac23 R_{nins}R_{nsnj} )x_n^4 +O(|x|^5).
\end{split}
\ee

Furthermore,
\begin{prop}[Proposition 3.2 of \cite{M1} or Lemma 2.4 of \cite{A}] \label{prop:marques}  At $P$, there hold
\begin{itemize}
\item[(1)] $h_{ij}=h_{ij,k}=h_{ij,kl}=0$,
\item[(2)] $R_{nn}=R_{nk}=\overline R_{kl}=0$,
\item[(3)] $R_{nn,n}=R_{nn,k}=R_{nk,l}=0$,
\item[(4)] $Sym_{klm}(\overline R_{kl,m})=0$,
\item[(5)] $R_{nn,nn}=-2(R_{ninj})^2$,
\item[(6)] $R_{nn,nk}=0$,
\item[(7)] $Sym_{kl}R_{nn,kl}=0$,
\item[(8)] $Sym_{klm}H_{, klm}=0$,
\item[(9)] $Sym_{klmp}(\frac{1}{2} \overline R_{kl,mp}+\frac{1}{9} \overline R_{rkls}\overline R_{rmps})=0,$
\item[(10)] $R=R_{,i}=R_{,n}=R_{,ni}=0,$
\item[(11)] $R_{,ii}=-\frac{1}{6}|\overline W|^2$
\item[(12)] $R_{ninj,ij}=-\frac 12 R_{,nn}-(R_{nknl})^2.$
\end{itemize}
Here $ \displaystyle{|\overline W|^2=\sum_{i,j,k,l=1}^{n-1}\overline W^2_{ijkl}}$ with $\overline W_{ijkl}$ being the components of the Weyl tensor $\overline W$ of $\pa M$ at $P$, and we also used a short notation that $\displaystyle{(R_{ninj})^2=\sum_{i,j=1}^{n-1}R^2_{ninj}}$.
\end{prop}

By the Proposition \ref{prop:marques}, and the second Bianchi identity, we have
\begin{align*}
\pa_i g^{ij}x_j&= R_{ninj,i}x_n^2x_j-\frac{1}{18}\overline R_{ikjs}\overline R_{imps} x_kx_mx_px_j +\frac 12 (R_{ninj,il}+R_{ninj,li})x_n^2x_lx_j \\
&\quad+\frac 16 \overline R_{jisl}R_{nsni}x_n^2x_lx_j +\frac 13 R_{ninj,ni}x_n^3x_j+O(|x|^5),\\
(g^{ij}-\delta_{ij}) x_ix_j&= R_{ninj}x_ix_jx_n^2+R_{ninj,k}x_n^2 x_ix_jx_k+\frac13 R_{ninj,n}x_n^3 x_ix_j
+\frac12 R_{ninj,kl}x_n^2x_ix_j x_k x_l\\&
\quad+\frac13 R_{ninj,nk} x_n^3x_ix_jx_k + (\frac{1}{12} R_{ninj,nn}+\frac23 R_{nins}R_{nsnj} )x_n^4x_ix_j +O(|x|^7),\\
(g^{ii}-\delta_{ii}) &= \frac{1}{18}\overline R_{ikls} \overline R_{imps}x_{k}x_l x_mx_p+\frac13 \overline R_{iksl} R_{nsni}x_n^2 x_k x_l+\frac{1}{2} ( R_{nins} )^2x_n^4 +O(|x|^5),\\
R&= \frac 12 R_{,nn} x_n^2 +\frac 12 R_{,ij}x_ix_j +O(|x|^3).
\end{align*}
Define
\begin{align*}
F[U_\lda](x)=(\pa_i g^{ij}\pa _j +(g^{ij}-\delta^{ij}) \pa_{ij} -\frac{n-2}{4(n-1)}R_g )U_\lda+O(|x|^{N-1}|\nabla U_\lda|).
\end{align*}
Recall
\[
\pa_j U_\lda= -(n-2)\lda^{\frac{n-2}{2}}\Big((\lda+x_n)^2+|x'|^2\Big)^{-\frac{n}{2}}x_j
\]
and
\[
\pa_{ij} U_{\lda}=  -(n-2)U_\lambda ((\lambda+x_n)^2+|x'|^2)^{-1}\delta_{ij}+n(n-2)U_\lambda ((\lambda+x_n)^2+|x'|^2)^{-2}x_ix_j.
\]
Then we have in $\Omega$ that
\be \label{eq:FU-2}
\begin{split}
&|F[U_\lda](x)| +|x+\lda e_n|^s|\nabla^s F[U_\lda](x)| 
\le  C\lda^{\frac{n-2}{2}} |x+\lda e_n|^{2-n},
\end{split}
\ee
where $s=1,2,3,4$, and $C>0$ depends only on $M,g,$ and $n$. Moreover, using symmetry,
\begin{align*}
&\frac{1}{|\mathbb{S}^{n-2}|}\int_{|x'|=r}F[U_\lda](x)\,\ud s_{x'}\\
&=\frac{n-2}{n-1} (\frac 12 R_{,nn}+(R_{ninj})^2)\frac{\lambda^{\frac{n-2}{2}}x_n^2r^2}{((\lambda+x_n)^2+r^2)^\frac{n}{2}}\\
&\quad +A+\frac{n(n-2)}{2(n-1)}(R_{ninj})^2\frac{\lda^{\frac{n-2}{2}}x_n^4r^2}{((\lambda+x_n)^2+r^2)^\frac{n+2}{2}}-\frac{n-2}{2} (R_{ninj})^2\frac{\lambda^{\frac{n-2}{2}}x_n^4}{((\lambda+x_n)^2+r^2)^\frac{n}{2}}\\
&\quad-\frac{n-2}{8(n-1)}R_{,nn}\frac{\lambda^{\frac{n-2}{2}}x_n^2}{((\lambda+x_n)^2+r^2)^\frac{n-2}{2}}+\frac{n-2}{48(n-1)^2}
|\overline W|^2\frac{\lambda^{\frac{n-2}{2}}r^2}{((\lambda+x_n)^2+r^2)^\frac{n-2}{2}}\\
&\quad\quad+O(|x|^3U_\lambda)\\
&=: \lambda^{\frac{n-2}{2}} g_\lambda(r,x_n)+O(|x|^3U_\lambda),
\end{align*}
where
\begin{align*}
A&=\frac{n(n-2)}{2|\mathbb{S}^{n-2}|}U_\lambda ((\lambda+x_n)^2+|x'|^2)^{-2}x_n^2\int_{|x'|=r} R_{ninj,kl}x_ix_j x_k x_l\,\ud s_{x'}\\
&=\frac{n(n-2)}{2|\mathbb{S}^{n-2}|}\frac{\lambda^{\frac{n-2}{2}}x_n^2}{((\lambda+x_n)^2+r^2)^\frac{n+2}{2}}\int_{|x'|=r}\sum_{i\neq j}(2R_{ninj,ij}+R_{nini,jj})x_i^2x_j^2+\sum_i R_{nini,ii}x_i^4\\
&=\frac{n(n-2)}{2|\mathbb{S}^{n-2}|}\frac{\lambda^{\frac{n-2}{2}}x_n^2}{((\lambda+x_n)^2+r^2)^\frac{n+2}{2}}\int_{|x'|=r}[\sum_{i\neq j}(\frac{2}{3}R_{ninj,ij}+\frac{1}{3}R_{nini,jj})+\sum_i R_{nini,ii}]x_i^4\\
&=\frac{n(n-2)}{(n+1)(n-1)}\frac{\lambda^{\frac{n-2}{2}}x_n^2r^4}{((\lambda+x_n)^2+r^2)^\frac{n+2}{2}}R_{ninj,ij}\\
&=-\frac{n(n-2)}{(n+1)(n-1)}[\frac 12 R_{,nn}+(R_{ninj})^2]\frac{\lambda^{\frac{n-2}{2}}x_n^2r^4}{((\lambda+x_n)^2+r^2)^\frac{n+2}{2}}.
\end{align*}
Here, we used \eqref{eq:trick} in the third and fourth equalities, and item (7) in Proposition \ref{prop:marques} in the fourth equality when we calculate $A$ in the above.

Let $W_\lda$ solve \eqref{Eq:reducedequation}. As in the previous section, we define
\[
W^{(1)}_{\lda}= \mathcal{G}(F[U_\lda]),  \quad  W^{(2)}_{\lda}= \mathcal{G}(F(W^{(1)}_{\lda})) \quad \mbox{and} \quad W_\lda^{(3)}= W_{\lda}-W^{(1)}_{\lda}-W^{(2)}_{\lda},
\]
where $\mathcal{G}(\cdot)$ is given in \eqref{eq:green-represent}.
By \eqref{eq:FU-2}, arguing as in the proof of Lemma \ref{lem:leading}, Lemma \ref{lem:leading2} and Lemma \ref{lem:leading3}, we have, for $n\ge 9$ and $s=1,2,3,4$,
\begin{equation}\label{eq:estimateofleading}
\begin{split}
|W_\lda^{(1)}|+|x+\lda e_n|^{s}|\nabla^s W_\lda^{(1)}|& \le  C \lda^{\frac{n-2}{2}} |x+\lda e_n|^{4-n},\\
|W_\lda^{(2)}|+|x+\lda e_n|^{s}|\nabla^s W_\lda^{(2)}|& \le  C \lda^{\frac{n-2}{2}} |x+\lda e_n|^{6-n},\\
|W_\lda^{(3)}|& \le  C \lda^{\frac{n-2}{2}} |x+\lda e_n|^{8-n}.
\end{split}
\end{equation}

We now proceed to show the second part of Theorem \ref{thm:main}.

\begin{proof}[Proof of Theorem \ref{thm:main} (ii)]
As in the previous section, we are going to calculate
\[
\int_{\Omega}U_\lda ^{\frac{n+2}{n-2}} W_\lda \,\ud x\quad\mbox{and}\quad \int_{\Omega}U_\lda ^{\frac{4}{n-2}} W_\lda^2 \,\ud x.
\]
Again,
\begin{align}
\int_{\Omega} U_{\lda}^{\frac{n+2}{n-2}} W_{\lda}&= \int_{\Omega} U_{\lda}^{\frac{n+2}{n-2}} W^{(1)}_{\lda}+ \int_{\Omega} U_{\lda}^{\frac{n+2}{n-2}} W^{(2)}_{\lda}+ \int_{\Omega} U_{\lda}^{\frac{n+2}{n-2}} W^{(3)}_{\lda} \nonumber \\
&= \int_{\Omega} \mathcal{G}( U_{\lda}^{\frac{n+2}{n-2}}) F[U_{\lda}]+\int_{\Omega} \mathcal{G}( U_{\lda}^{\frac{n+2}{n-2}}) F[W^{(1)}_{\lda}]+O(\lda^{6})\label{eq:decomp2}.
\end{align}

\emph{Step 1:} We first calculate 
\[
\int_{\Omega} \mathcal{G}( U_{\lda}^{\frac{n+2}{n-2}}) F[U_{\lda}].
\]

Using symmetry, we have
\begin{align}
\int_{\Omega} \mathcal{G}( U_{\lda}^{\frac{n+2}{n-2}}) F[U_{\lda}]&
= \lda^{\frac{n-2}{2}}\int_{\Omega} \mathcal{G}( U_{\lda}^{\frac{n+2}{n-2}}) g_{\lda}(|x'|,x_n)\,\ud x +O(\lda^{5}) \nonumber\\&
=\lda^4\int_{\Omega_\lda}  \mathcal{G}_\lda( U_{1}^{\frac{n+2}{n-2}}) g_{1}(|x'|,x_n)\,\ud x +O(\lda^{5}) \nonumber\\&
=\frac{1}{2n}\lda^4\int_{\R^n_+}  y_n((1+y_n)^2+|y'|^2)^{-\frac{n}{2}} g_1(|y'|,y_n)\,\ud y+O(\lda^{n-1})+O(\lda^{5}) \nonumber\\&
=:\frac{1}{2n}\lda^4\sum_{j=1}^6 I_6 +O(\lda^{5}),
\label{eq:umbilical-a}
\end{align}
where we have used \eqref{eq:theleadingterm} in the third equality, and
\begin{align*}
I_1&=\frac{n-2}{n-1} (\frac 12 R_{,nn}+(R_{ninj})^2)\int_{\R^n_+}  \frac{y_n^3|y'|^2}{((1+y_n)^2+|y'|^2)^n},\\
I_2&=-\frac{n(n-2)}{(n+1)(n-1)}[\frac 12 R_{,nn}+(R_{ninj})^2]\int_{\R^n_+} \frac{y_n^3|y'|^4}{((1+y_n)^2+|y'|^2)^{n+1}},\\
I_3&=\frac{n(n-2)}{2(n-1)}(R_{ninj})^2\int_{\R^n_+}\frac{y_n^5|y'|^2}{((1+y_n)^2+|y'|^2)^{n+1}},\\
I_4&=-\frac{n-2}{2} (R_{ninj})^2\int_{\R^n_+}\frac{y_n^5}{((1+y_n)^2+|y'|^2)^n},\\
I_5&=-\frac{n-2}{8(n-1)}R_{,nn}\int_{\R^n_+}\frac{y_n^3}{((1+y_n)^2+|y'|^2)^{n-1}},\\
I_6&=\frac{n-2}{48(n-1)^2}|\overline W|^2\int_{\R^n_+}\frac{y_n|y'|^2}{((1+y_n)^2+|y'|^2)^{n-1}}.
\end{align*}
To calculate $I_1-I_6$, we will need the following identities, which can be obtained by integration by parts or change of variables.
\begin{align*}
\int_0^\infty\frac{r^{n+2}}{(1+r^2)^{n+1}}\,\ud r&=\frac{n+1}{2n}\int_0^\infty\frac{r^n}{(1+r^2)^n}\,\ud r,\\
\int_0^\infty\frac{r^{n-2}}{(1+r^2)^{n}}\,\ud r&=\int_0^\infty\frac{r^n}{(1+r^2)^n}\,\ud r,\\
\int_0^\infty\frac{r^{n-2}}{(1+r^2)^{n-1}}\,\ud r&=2\int_0^\infty\frac{r^n}{(1+r^2)^n}\,\ud r,\\
\int_0^\infty\frac{r^n}{(1+r^2)^{n+1}}\,\ud r&=\frac{n-1}{2n}\int_0^\infty\frac{r^n}{(1+r^2)^n}\,\ud r.
\end{align*}
By calculations, we have
\begin{align*}
I_1&=\frac{n-2}{n-1} (\frac 12 R_{,nn}+(R_{ninj})^2)\frac{6|\mathbb S^{n-2}|}{(n-2)(n-3)(n-4)(n-5)}\int_0^\infty\frac{r^n}{(1+r^2)^n}\,\ud r,\\
I_2&=-\frac{n(n-2)}{(n+1)(n-1)}[\frac 12 R_{,nn}+(R_{ninl})^2] \frac{3(n+1)|\mathbb S^{n-2}|}{n(n-2)(n-3)(n-4)(n-5)}\int_0^\infty\frac{r^n}{(1+r^2)^n}\,\ud r,\\
I_3&=\frac{n(n-2)}{2(n-1)}(R_{ninj})^2\frac{120|\mathbb S^{n-2}|}{n(n-1)(n-2)(n-3)(n-4)(n-5)}\frac{n-1}{2n}\int_0^\infty\frac{r^n}{(1+r^2)^n}\,\ud r,\\
I_4&=-\frac{n-2}{2} (R_{ninj})^2\frac{120|\mathbb S^{n-2}|}{n(n-1)(n-2)(n-3)(n-4)(n-5)}\int_0^\infty\frac{r^n}{(1+r^2)^{n}}\,\ud r,\\
I_5&=-\frac{n-2}{8(n-1)}R_{,nn}\frac{12|\mathbb S^{n-2}|}{(n-2)(n-3)(n-4)(n-5)}\int_0^\infty\frac{r^n}{(1+r^2)^n}\,\ud r,\\
I_6&=\frac{n-2}{48(n-1)^2}|\overline W|^2\frac{|\mathbb S^{n-2}|}{(n-4)(n-5)}\int_0^\infty\frac{r^n}{(1+r^2)^{n-1}}\,\ud r.
\end{align*}
Suppose that  $$|\overline W|^2+(R_{ninj})^2\neq 0,$$ which is equivalent to the full Weyl tensor $$|W|\neq 0$$ at $P$ as explained in \cite{M1}. We have
\be  \label{eq:um-A}
\begin{split}
\sum_{i=1}^6 I_i&=\frac{3(n-10)(R_{ninj})^2|\mathbb S^{n-2}|}{n(n-1)(n-3)(n-4)(n-5)}\int_0^\infty\frac{r^n}{(1+r^2)^{n}}\,\ud r\\
&\quad +\frac{n-2}{48(n-1)^2}|\overline W|^2\frac{|\mathbb S^{n-2}|}{(n-4)(n-5)}\int_0^\infty\frac{r^n}{(1+r^2)^{n-1}}\,\ud r,
\end{split}
\ee
which is nonnegative if $n\ge 10$, and positive  if $n\ge 11$.

\emph{Step 2:} Next, we estimate
\[
\int_{\Omega} \mathcal{G}( U_{\lda}^{\frac{n+2}{n-2}}) F[W^{(1)}_{\lda}].
\]

It follows from \eqref{eq:estimateofleading} and Proposition \ref{prop:marques} that
\begin{align*}
F[W^{(1)}_{\lda}]=\frac{1}{3}\overline{R}_{ikjl}x_kx_l\pa_{ij}W^{(1)}_{\lda}+R_{ninj} x_n^2\pa_{ij}W^{(1)}_{\lda}+O(\lda^{\frac{n-2}{2}}|x+\lda e_n|^{5-n}).
\end{align*}
Also, using item (2) in Proposition \ref{prop:marques} and the fact that $U_\lda$ is radial in $x'$, we have
\begin{align*}
W^{(1)}_\lda= \mathcal{G} (F[U_\lda])&= \mathcal{G} \Big(R_{ninj}  x_n^2 \pa_{ij}U_{\lda} +O(|x|U_\lda)\Big)\\
&=\mathcal{G} \Big(R_{ninj}  x_n^2 \pa_{ij}U_{\lda}\Big) +O(\lda^{\frac{n-2}{2}}|x+\lda e_n|^{5-n}).
\end{align*}


Hence, using \eqref{eq:theleadingtermscaled}, \eqref{eq:estimateofleading}, integration by parts and $R_{nn}=0$, we have
\begin{align*}
&\int_{\Omega}  \mathcal{G}( U_{\lda}^{\frac{n+2}{n-2}}) R_{ninj} x_n^2\pa_{ij}W^{(1)}_\lda\\
&=\frac{1}{2n}\lda^{n/2}\int_{\Omega} x_n^3|x+\lda e_n|^{-n}  R_{ninj} \pa_{ij}W^{(1)}_\lda+O(\lda^{n-1})\\
&=\frac{n+2}{2}\lda^{n/2}\int_{\Omega}x_n^3|x+\lda e_n|^{-n-4}R_{ninj}x_ix_j \mathcal{G}(R_{nknl} y_n^2  \pa_{kl}U_{\lda})+O(\lda^5)\\
&=\frac{n+2}{2}\lda^{4}\int_{\Omega_\lda}x_n^3|x+ e_n|^{-n-4}R_{ninj}x_ix_j \mathcal{G}_\lda(R_{nknl} y_n^2  \pa_{kl}U_{1})+O(\lda^5),
\end{align*}
where the integrals on the boundary $\pa\Omega$ coming from the integration by parts are absorbed in $O(\lda^5)$. Similarly, 
\begin{align*}
&\int_{\Omega}  \mathcal{G}(U_{\lda}^{\frac{n+2}{n-2}}) \overline{R}_{ikjl}x_kx_l\pa_{ij}W^{(1)}_\lda\\
&=\frac{1}{2n}\lda^{n/2}\int_{\Omega} x_n|x+\lda e_n|^{-n}  \overline{R}_{ikjl}x_kx_l \pa_{ij}W^{(1)}_\lda+O(\lda^{n-1})\\
&=\frac{1}{2n}\lda^{n/2}\int_{\Omega} x_n W^{(1)}_\lda\pa_{ij}\Big( |x+\lda e_n|^{-n}  \overline{R}_{ikjl}x_kx_l\Big)+O(\lda^{n-1})\\
&=O(\lda^{n-1}),
\end{align*}
where item (2) in Proposition \ref{prop:marques} is used in the last equality. 

Therefore,
\begin{align*}
&\int_{\Omega}  \mathcal{G}( U_{\lda}^{\frac{n+2}{n-2}}) F[W^{(1)}_{\lda}]\\
&=\int_{\Omega}  \mathcal{G}( U_{\lda}^{\frac{n+2}{n-2}}) \Big(R_{ninj} x_n^2\pa_{ij}W^{(1)}_\lda+O(\lda^{\frac{n-2}{2}}|x+\lda e_n|^{5-n})\Big)+O(\lda^{n-1})\\
& =\frac{n+2}{2}\lda^{4}\int_{\Omega_\lda}x_n^3|x+ e_n|^{-n-4}R_{ninj}x_ix_j \mathcal{G}_\lda(R_{nknl} y_n^2  \pa_{kl}U_{1})+O(\lambda^5).
\end{align*}

We are going to estimate  $\mathcal{G}_\lda(R_{nknl} y_n^2  \pa_{kl}U_{1})(x).$
Recall the notation that
\[
\overline{\mathcal{G}} (f)=\int_{\R^n_+}\overline G(x,y)f(y)\,\ud y.
\]
Then
\begin{equation}\label{eq:auxupperbound2}
|\overline{\mathcal{G}}( y_n^2\pa_{ij} U_{1} )(x)|\le C\int_{\R^n_+}\frac{1}{|x-y|^{n-2}|y+ e_n|^{n-2}}\,\ud y\le C|x+e_n|^{4-n}.
\end{equation}
By the maximum principle for harmonic function, we have
\begin{equation}\label{eq:aux222}
\|\mathcal{G}_\lambda (y_n^2\pa_{ij} U_{1} )-\overline{\mathcal{G}}(y_n^2\pa_{ij} U_{1} )\|_{L^\infty(\Omega_\lambda)}\le \|\overline{\mathcal{G}}(y_n^2\pa_{ij} U_{1} )\|_{L^\infty(\partial\Omega_\lambda)}\le C\lambda^{n-4}.
\end{equation}

Hence,
\begin{align*}
&\int_{\Omega}  \mathcal{G}( U_{\lda}^{\frac{n+2}{n-2}}) F[W^{(1)}_{\lda}]\\
&=\frac{n+2}{2}\lda^{4}\int_{\Omega_\lda}x_n^3|x+ e_n|^{-n-4}R_{ninj}x_ix_j \overline{\mathcal{G}}(R_{nknl} y_n^2  \pa_{kl}U_{1})+O(\lda^5)\\
&=\frac{n+2}{2}\lda^{4}\int_{\R^n_+}x_n^3|x+ e_n|^{-n-4}R_{ninj}x_ix_j \overline{\mathcal{G}}(R_{nknl} y_n^2  \pa_{kl}U_{1})+O(\lda^5)\\
&=\frac{n^3-4n}{2}\lambda^{4}\int_{\R^n_+} |x+e_n|^{-n-4} R_{ninj}x_ix_j x_n^3\overline{\mathcal{G}}(R_{nknl}y_ky_ly_n^2|y+ e_n|^{-2-n})+O(\lambda^5).
\end{align*}

Now we calculate $\Psi:= \overline{\mathcal{G}}(R_{nknl}y_ky_ly_n^2|y+ e_n|^{-2-n})$.
It solves
\[
-\Delta \Psi= R_{nknl} x_kx_l \xi(|x'|,x_n) \quad \mbox{in }\R^n_+, \quad \Psi= 0  \quad \mbox{on }\pa \R^n_+,
\]
where $\xi(|x'|,x_n)= x_n^2  (|x'|^2+ (x_n+1)^2)^{-\frac{2+n}{2}}$. We look for the unique  $\Psi= R_{nknl}x_kx_l  \Lda(|x'|,x_n)$ for some $\Lda=\Lda(r,s)$. Then $\Lda$ satisfies
\begin{equation}\label{eq:highlap-12}
-\partial_{rr} \Lda-\frac{n+2}{r}\partial_r\Lda-\partial_{ss}\Lda=\xi(r,s)  \quad \Lda(r,0)= 0\quad \mbox{for }r,s>0.
\end{equation}
where $R_{nn}=0$ is used.  As before, the operator on the left hand side of \eqref{eq:highlap-0} can be considered as the Laplacian operator $-\Delta$ in $\R^{n+4}$ applying to $\Lda$ which is radial in the first $n+3$ coordinates. Therefore, for $\tilde \Lda(z)=\Lda(|z'|,z_{n+4})$ with $z'=(z_1,\cdots,z_{n+3})$, it satisfies
\begin{equation}\label{eq:highlap-13}
-\Delta \tilde \Lda(z)= z^2_{n+4} |z+ e_{n+4}|^{-2-n}  \quad \mbox{in }\R^{n+4}_+, \quad \tilde \Lda= 0  \quad \mbox{on }\pa\R^{n+4}_+.
\end{equation}
This $\tilde \Lda$ can be solved using Green's function of $-\Delta$ in $\R^{n+4}_+$, and thus, $\tilde \Lda>0$ everywhere.  Moreover, by the same argument as of \eqref{eq:solveleading}, we have that
\[
\tilde \Lda(z)\le C   |z+ e_{n+4}|^{2-n} \quad \mbox{in }\R^{n+4}_+,
\]
and  for all $r,s\ge 0$,
\[
0<\Lda(r,s)\le C  (r^2+(1+s)^2)^\frac{2-n}{2}.
\]

Hence,
\[
\int_{\Omega}  \mathcal{G}( U_{\lda}^{\frac{n+2}{n-2}}) F[W^{(1)}_{\lda}]=\frac{n^3-4n}{2}\lambda^{4}\int_{\R^n_+} |x+e_n|^{-(n+4)} x_n^3 \Lda(|x'|,x_n) (R_{ninj}x_ix_j)^2 +O(\lambda^5)
\]
Since $\sum_{i} R_{nini}=R_{nn}=0$, we have
\begin{align*}
\int_{\mathbb S^{n-2}}(R_{ninj}x_ix_j)^2 &= \int_{\mathbb S^{n-2}} \sum_{i}R_{nini}^2x_i^4+2\sum_{i\neq j}R_{ninj}^2x_i^2x_j^2+\sum_{i\neq j}R_{nini}R_{njnj}x_i^2x_j^2\\
&=\int_{\mathbb S^{n-2}} x_1^4\Big(\sum_{i}R_{nini}^2+\frac{2}{3}\sum_{i\neq j}R_{ninj}^2+\frac 13\sum_{i\neq j}R_{nini}R_{njnj}\Big)\\
&=\frac{2}{3} (R_{ninj})^2\int_{\mathbb S^{n-2}} x_1^4,
\end{align*}
where in the second equality we used \eqref{eq:trick}.

Hence,
\begin{align}
\int_{\Omega}  \mathcal{G}( U_{\lda}^{\frac{n+2}{n-2}}) F[W^{(1)}_{\lda}]&=\frac{n^3-4n}{3}\lambda^{4}(R_{ninj})^2\int_{\R^n_+} |x+e_n|^{-(n+4)} x_n^3 \Lda(|x'|,x_n) x_1^4 +O(\lambda^5)\nonumber\\
&=C(n)(R_{ninj})^2\lda^4+O(\lambda^5)\label{eq:estimateofleadingterm22}
\end{align}
for some $C(n)>0$.

\emph{Step 3:} Finally, we estimate
\[
\int_{\Omega}U_\lda ^{\frac{4}{n-2}} W_\lda^2 \,\ud x.
\]

It follows from \eqref{eq:estimateofleading} that
\[
\int_{\Omega}U_\lda ^{\frac{4}{n-2}} W_\lda^2 \,\ud x=\int_{\Omega}U_\lda ^{\frac{4}{n-2}} (W^{(1)}_\lda)^2 \,\ud x+O(\lda^5).
\]
By the definition of $W^{(1)}_\lda$, for $n\ge 10$ we have
\begin{align*}
W^{(1)}_\lda(x)&=  \mathcal{G} \Big(  R_{ninj}  y_n^2 \pa_{ij}U_{\lda} +O(|y|U_\lda)\Big)(x)\\
&=\mathcal{G}(  R_{ninj}  y_n^2 \pa_{ij}U_{\lda} )(x) +O(\lda^{\frac{n-2}{2}}|x+\lda e_n|^{5-n}).
\end{align*}

Therefore, using \eqref{eq:auxupperbound2} and \eqref{eq:aux222}, we have
\begin{align}
\int_{\Omega} U_\lda ^{\frac{4}{n-2}} (W^{(1)}_\lda)^2&= \lambda^4 \int_{\Omega_\lambda} U_1 ^{\frac{4}{n-2}} (\mathcal{G}_\lda (R_{ninj}   y_n^2 \pa_{ij}U_{1}))^2 + O(\lda^5)\nonumber\\
&= \lambda^4 \int_{\Omega_\lambda} U_1 ^{\frac{4}{n-2}} (  \overline{\mathcal{G}} (R_{ninj} y_n^2 \pa_{ij}U_{1}))^2 + O(\lda^5)\nonumber\\
&= \lambda^4 \int_{\R^n_+} U_1 ^{\frac{4}{n-2}} (\overline{\mathcal{G}} (R_{ninj} y_n^2 \pa_{ij}U_{1}))^2 + O(\lda^5)\nonumber\\
&= n^2(n-2)^2\lambda^4 \int_{\R^n_+} U_1 ^{\frac{4}{n-2}} (\overline{\mathcal{G}}(R_{ninj} y_iy_jy_n^2|y+e_n|^{-2-n}))^2  + O(\lda^5) \nonumber \\
&=n^2(n-2)^2\lambda^4 \int_{\R^n_+} U_1 ^{\frac{4}{n-2}} \Lda^2(|x'|,x_n) (R_{ninj}   x_ix_j)^2 + O(\lda^5) \nonumber \\
&=\frac{2n^2(n-2)^2}{3}\lambda^4 (R_{ninj})^2\int_{\R^n_+} U_1 ^{\frac{4}{n-2}}\Lda^2(|x'|,x_n) x_1^4 +O(\lda^5)  \nonumber \\
&=C(n)\lambda^4 (R_{ninj})^2 +O(\lda^5), \label{eq:28}
\end{align}
for some $C(n)>0$, where $\Lda$ is the same one as in \eqref{eq:highlap-12}.

It follows that, using \eqref{eq:decomp2}, \eqref{eq:umbilical-a}, \eqref{eq:um-A}, \eqref{eq:estimateofleadingterm22} and \eqref{eq:28},
\begin{align*}
&\int_{M} |u_\lda|^{\frac{2n}{n-2}}\,\ud vol_g \nonumber \\& =  \int_{\Omega} (U_\lda+W_\lda)^{\frac{2n}{n-2}}(1+O(|x|^N))\,\ud x +O(\lda^n)  \\&
=   \int_{\Omega} \left( U_\lda ^{\frac{2n}{n-2}}  +\frac{2n}{n-2}U_\lda ^{\frac{n+2}{n-2}} W_\lda +\frac{n(n+2)}{(n-2)^2}U_\lda ^{\frac{4}{n-2}} W_\lda^2 +O(|W_\lda|^{\frac{2n}{n-2}})\right)\,\ud x  +O(\lda^5)\\
&=   \int_{\Omega} \left( U_\lda ^{\frac{2n}{n-2}}  +\frac{2n}{n-2}U_\lda ^{\frac{n+2}{n-2}} W_\lda +\frac{n(n+2)}{(n-2)^2}U_\lda ^{\frac{4}{n-2}} (W^{(1)}_\lda)^2 \right)\,\ud x  +O(\lda^{\frac{4n}{n-2}})\\
&= 2^{-n} \w_n +C(n)\lda^4 +O(\lda^{\frac{4n}{n-2}}),
\end{align*}
where $ C(n)>0$ provided that $n\ge 10$ and $|W|^2>0$.
Using \eqref{eq:denominator-part}, we complete the proof of Theorem \ref{thm:main}.
\end{proof}

\begin{rem}\label{rem:2}
From Steps 2 and 3 in this section, if one can explicitly calculate $\Lda$ in \eqref{eq:highlap-12}, or equivalently $\tilde \Lda$ in \eqref{eq:highlap-13}, then the dimension assumption might be reduced.
\end{rem}

\begin{proof}[Proof of Corollary \ref{cor:1}] We only need to prove the ``only if" part. Since $n\ge 12$, it follows from Theorem \ref{thm:main} that $\pa \mathcal{O}$ is umbilic. Since $\pa \mathcal{O}$ is smooth and connected, $\pa \mathcal{O}$ has to be a sphere. 
\end{proof}

\bigskip

\noindent T. Jin

\noindent Department of Mathematics, The Hong Kong University of Science and Technology\\
Clear Water Bay, Kowloon, Hong Kong\\[1mm]
Email: \textsf{tianlingjin@ust.hk}

\medskip

\noindent J. Xiong

\noindent School of Mathematical Sciences, Beijing Normal University\\
Beijing 100875, China\\[1mm]
Email: \textsf{jx@bnu.edu.cn}


\small

\begin{thebibliography}{99}

\bibitem{A} S.  Almaraz,
        \emph{An existence theorem of conformal scalar-flat metrics on manifolds with boundary.}
        Pacific J. Math. \textbf{248} (2010), no. 1, 1--22.

\bibitem{Aubin} T. Aubin,
        \emph{\'Equations différentielles non lin\'eaires et probl\`eme de Yamabe concernant la courbure scalaire}. J. Math. Pures Appl. (9) \textbf{55} (1976), no. 3, 269--296.

\bibitem{Carleman} T. Carleman,
        \emph{Zur Theorie der Minimalfl\"achen.}
        Math. Z. \textbf{9} (1921),  154--160.

\bibitem{E} J.F.  Escobar,
        \emph{Conformal deformation of a Riemannian metric to a scalar flat metric with constant mean curvature on the boundary.}
        Ann. of Math. (2) \textbf{136} (1992), no. 1, 1--50.

\bibitem{GZ} M. Gluck and M. Zhu,
             \emph{An extension operator on bounded domains and applications.} Preprint.

\bibitem{HWY1} F. Hang, X. Wang and X. Yan,
         \emph{Sharp integral inequalities for harmonic functions.}
         Comm. Pure Appl. Math.  \textbf{61} (2008), no. 1, 54--95.

\bibitem{HWY}  F. Hang, X. Wang and X. Yan,
         \emph{An integral equation in conformal geometry.}
         Ann. Inst. H. Poincar\'e Anal. Non Lin\'eaire \textbf{26} (2009), no. 1, 1--21.

\bibitem{J} S. Jacobs,
          \emph{An isoperimetric inequality for functions analytic in multiply connected domains.}
          Mittag-Leffler Institute report, 1972.


\bibitem{M1} F. Marques,
         \emph{Existence results for the Yamabe problem on manifolds with boundary.}
         Indiana Univ. Math. J. \textbf{54} (2005), no. 6, 1599--1620.

\bibitem{M2} F. Marques,
         \emph{Conformal deformations to scalar-flat metrics with constant mean curvature on the boundary.}
         Comm. Anal. Geom. \textbf{15} (2007), no. 2, 381--405.

\bibitem{Schoen} R. Schoen,
         \emph{Conformal deformation of a Riemannian metric to constant scalar curvature}. J. Differential Geom. \textbf{20} (1984), no. 2, 479--495.

\bibitem{Trudger} N.S. Trudinger, \emph{Remarks concerning the conformal deformation of Riemannian structures on compact manifolds}. Ann. Scuola Norm. Sup. Pisa (3) \textbf{22} (1968),  265--274.

\bibitem{X2} J. Xiong,
         \emph{On a conformally invariant integral equation involving Poisson kernel}. Preprint.

\bibitem{Yamabe} H. Yamabe,
         \emph{On a deformation of Riemannian structures on compact manifolds}. Osaka Math. J. \textbf{12} (1960), 21--37.
\end{thebibliography}
\end{document}